\documentclass[11pt]{article}

\usepackage{amsmath,amssymb}

\newcommand{\comment}[1]{}

\newtheorem{thm}[subparagraph]{Theorem}
\newtheorem{lem}[subparagraph]{Lemma}

\newtheorem{cor}[subparagraph]{Corollary}
\newcommand{\Pow}{\textsf{Pow}}

\newcommand{\be}{\begin}
\newcommand{\en}{\end}

\begin{document}
\title{On the existence of Stone-\v{C}ech compactification}

\author{Giovanni Curi \\ \\
{\small Dipartimento di Informatica - Universit\`a di Verona}\\
{ \small  Strada le Grazie 15 - 37134 Verona, Italy.}\\
{ \small e-mail: giovanni.curi@univr.it}\\ \\ }
\date{\small January 9, 2010.}

\maketitle

\paragraph* {Introduction.} In 1937 E. \v{C}ech
and M.H. Stone, independently,
introduced the maximal compactification  of a completely regular topological space,
thereafter called Stone-\v{C}ech compactification \cite{Cech,Stone}.
In the introduction of \cite{Cech} the non-constructive character of this result is so described:
``It must be emphasized that $\beta(S)$ [the Stone-\v{C}ech compactification of $S$] may be defined only formally (not constructively) since it exists only in virtue of Zermelo's theorem''.\footnote{Over ZF, the existence of Stone-\v{C}ech compactification has later been proved  equivalent to the Prime Ideal Theorem \cite{BM80}.}

By replacing topological spaces with locales, Banaschewski and Mulvey \cite{BM80,BM84,BM}, and Johnstone  \cite{J82} obtained choice-free intuitionistic proofs of  Stone-\v{C}ech compactification.
Although valid in any topos, these localic constructions  rely - essentially, as is to be demonstrated -  on highly impredicative principles, and thus cannot be considered as constructive in the sense of the main systems for constructive mathematics, such as Martin-L\"of's constructive type theory  and Aczel's constructive set theory.

In \cite{CuSC} I  characterized the locales of which the Stone-\v{C}ech compactification can be defined in constructive type theory CTT, and in the formal system CZF+uREA+DC, a  natural extension of Aczel's system for constructive set theory CZF by a strengthening of the Regular Extension Axiom REA and the principle of Dependent Choice.

In this paper I show that this characterization continues to hold over the standard system CZF$^+$, i.e.,  CZF plus REA,  thus removing in particular any dependency from a choice principle.
This will follow by a result of independent interest, namely
the proof that the class of continuous mappings from a compact regular locale $X$ to a regular a set-presented locale $Y$ is a set in CZF (even without REA).

By exploiting the consistency of CZF+REA with a general form of Troelstra's \emph{principle of uniformity} \cite{TvD}, I then prove that the existence of Stone-\v{C}ech compactification of a non-degenerate  Boolean locale is independent of the axioms of CZF (+REA), so that the aforementioned  characterization  characterizes a proper subcollection  of the collection of all locales.
The same also holds for several, even impredicative, extensions of CZF+REA, as well as for CTT.
This  is in strong contrast with what happens in the context of Higher-order Heyting arithmetic HHA - and thus in any topos-theoretic universe: within HHA  the constructions in  \cite{BM80,BM84,BM, J82} of Stone-\v{C}ech compactification  can be carried out for every locale.

\medskip

\paragraph {Constructive Set Theory and the Principle of Uniformity.}\label{Homsets} \setcounter{subparagraph}{0} The reader is invited to consult \cite{AR,ML84} for background on Aczel's constructive set theory, CZF, and Martin-L\"of's type theory, CTT. Here I shall be working in constructive set theory, but the arguments to be presented may be adapted also to the setting of CTT.

The system CZF is a subsystem of ZF that uses  intuitionistic rather than classical logic, has only  a restricted form of the Separation Scheme,  and  does not have the Powerset Axiom.
More specifically, CZF is based on first-order intuitionistic logic with equality, has  $\in$ as the only non-logical symbol,  and  has Extensionality, Pair, Union, Infinity, Set Induction,  Restricted Separation (i.e., Separation for bounded formulae), Strong Collection and Subset Collection as non-logical axioms and schemes.
Strong Collection is the following scheme:

\be{description}
\item[Strong Collection] For every set $a$,  if $(\forall x\in
a)(\exists y)\; \phi(x,y)$, then there is a set $b$ such that $(\forall x\in a)(\exists y\in b)\; \phi(x,y)$ and
$(\forall y\in b)(\exists x\in a)\; \phi(x,y)$.
\en{description}

\noindent A purely logical consequence of Strong Collection, to be often applied in the following (sometimes tacitly), is the usual Replacement Scheme.

\be{description}
\item[Replacement] For every set $a$,  if $(\forall x\in
a)(\exists ! y)\; \phi(x,y)$, then there is a set $b$ such that $\forall y(y\in b\leftrightarrow (\exists x\in a)\; \phi(x,y))$.
\en{description}

The Subset Collection Scheme is a strengthening of Myhill's Exponentiation Axiom, asserting that the class of functions between any two sets is a set.  Subset Collection is often presented in the equivalent form (over the remaining axioms of CZF)  of the Fullness Axiom.
For sets $a,b,$ let ${\bf mv}(b^a)$ be the class of  subsets $r$ of $a\times b$ such that $(\forall x\in a)(\exists
y\in b)\; (x,y)\in r$.

\be{description}
\item[Fullness] Given sets $a,b$ there is a subset $c$ of ${\bf
  mv}(a^b)$ such that for every  $r\in {\bf mv}(a^b)$ there is $r_0\in  c$ with $r_0\subseteq r$.
\en{description}

\noindent This principle may in some cases be used to replace applications of the fully impredicative Powerset Axiom. The system one obtains  replacing Subset Collection with the Exponentiation Axiom is usually denoted by CZF$_{exp}$.

The Regular Extension Axiom REA, stating that every set is a subset of a regular set, is often added to the axioms of CZF;  the resulting system is denoted by CZF$^+$. REA is needed in order to prove that certain inductively defined classes are sets (see \cite{AR} for more information).

Both CZF+REA and CTT can consistently be extended with a general form of Troelstra's principle of uniformity (van den Berg \& Moerdijk \cite{BergMoerdijk}, Coquand, cf. \cite{CuriPA}). In constructive set theory, this is formulated as the following schema:

\be{description}
\item[GUP] For every set $a$, if $(\forall x)(\exists y\in a)\phi(x,y)$, then $(\exists y\in a)(\forall x)\phi(x,y).$
\en{description}

\noindent This  principle may be seen as resulting  from the particular case in which $a$ is the set $\omega$  of natural numbers, i.e., from  the uniformity principle in its standard form, and the principle that every set is subcountable, also consistent with CZF+REA  (\cite{BergMoerdijk}; see also \cite{Lubarsky2006,Rathjen94,Rathjen02,Rathjen06,Streicher}).
In fact, several extensions of CZF, as CZF+REA+PA+Sep, where PA is the Presentation Axiom (implying the Dependent and Countable Choice principles), and Sep is the impredicative full Separation Scheme, have been proved  consistent with GUP  \cite{BergMoerdijk}.

As is customary in classical set theory, class notation and terminology can be exploited in this context  \cite{AR}. Recall also that a set is \emph{finite} in this setting if it (is empty, or) can be finitely enumerated (possibly with repetitions), and that  the class $\Pow_{fin}(S)$ of  finite subsets of a set $S$ is a set in CZF (and CZF$_{exp}$).

We conclude this introductory section by proving a first important general consequence of the consistency of the systems we are considering with GUP.
A  (large)  \emph{$\bigvee$-semilattice} is a partially ordered class (i.e., a class together with a class relation on it satisfying the usual axioms for a partial order) that  has  suprema for arbitrary \emph{subsets}.  A $\bigvee$-semilattice need not be a (large) \emph{complete lattice}, i.e., need not have also infima of arbitrary subsets. A  \emph{class-frame}, or \emph{class-locale}, $X$ is a $\bigvee$-semilattice that has a top element $\top$,  binary meets, and that is such that meets  distribute over suprema of arbitrary sets of elements of $X$ \cite{A}. A \emph{class-preframe} is defined in the same way as a class-frame, but suprema have to exist only of directed subsets, and meets  are required  to distribute only over these suprema. Note that a preframe need not have a smallest element. A partially ordered class is \emph{degenerate} if it consists of a single element. Then, a $\bigvee$-semilattice  $L$ is degenerate iff  $L=\{\bot\}$, with $\bot=\bigvee \emptyset$, while a class-preframe $P$ is degenerate iff $P=\{\top\}$.


Lattices of the above kinds arise everywhere in mathematics; the fact that they are carried by sets is often an essential tacit assumption in the theory and applications of these structures, e.g. when  classes of ideals are considered.
A consequence of the consistency of CZF  with the generalized uniformity principle is that, constructively, in no non-trivial case this assumption is legitimate.

\begin{lem}[The Main Lemma]\label{Noframeset}
No non-degenerate $\bigvee$-semilattice,  or class-preframe, and hence no
non-degenerate class-frame, can be proved to have a set of elements in
CZF (+REA+PA+Sep).
\end{lem}

\noindent \textbf{Proof.} Let $L$ be a non-degenerate $\bigvee$-semilattice, and assume $L$ is carried by a set. Then the class $\{x\in L: \emptyset\in y\}$ is a set for every set $y$. Therefore, $(\forall y)(\exists a\in L) a= \bigvee\{x\in L: \emptyset\in y\}$. In  CZF (+REA+PA+Sep)+GUP one then gets $(\exists a\in L)(\forall y) a= \bigvee\{x\in L: \emptyset\in y\}$, so that $L$ must be degenerate, as follows by  first taking $y=\emptyset$, then $y=\{\emptyset\}$. So $L$ is not a set in
CZF (+REA+PA+Sep)+GUP, and thus cannot be proved to be a set in CZF (+REA+PA+Sep).

The proof for $P$ a non-degenerate preframe (and in fact for any non-degenerate partially ordered class with a greatest element and joins of directed subsets) is similar, but one considers instead the set $\{x\in P: \emptyset\in y\}\cup \{b\}$, for $b\in P$. It is an easy exercise in intuitionistic logic to check that, for every $y$, this set is directed, so that it has a join in $P$ for every $y$. Reasoning as in the previous case, one gets that $b=\top$. As this holds for every $b\in P$, it follows that $P$ is degenerate, against the hypothesis.


\medskip

\noindent
Clearly, given any set $X$, its powerclass $\Pow(X)$ with  intersection as meet and union of arbitrary set-indexed families of subsets as join is a
frame. Therefore, no non-trivial instance of the Powerset Axiom is constructively derivable.

\begin{cor}
For no non-empty set $X$ the powerclass $\Pow(X)$ can be proved to form a set in CZF(+REA+PA+Sep).
\end{cor}

\paragraph {Small homsets of continuous maps in CZF.}\setcounter{subparagraph}{0}
Locales, or formal spaces, provide a suitable substitute to the concept of topological space in choice-free and/or intuitionistic settings \cite{J82,FS79}. In CZF, or CTT, due to the absence of powersets,  the concept of locale needs to be formulated with special care \cite{A,CSSV} (recall also the Main Lemma). In CZF, a class-locale is said to be \emph{set-generated} by a subclass $B$ if:
\begin{itemize}

\item [$i.$] $B$ is a set,

\item [$ii.$] the class $\{b\in B:b\leq x\}$ is a set and $x=\bigvee\{b\in B:b\leq
x\}$, for all  $x\in X$.

\end{itemize}

In a fully impredicative context  as intuitionistic set theory IZF,  set-generated class-locales and ordinary locales come to the same thing. Here a set-generated class-locale $(X,B)$ will simply be referred to as a \emph{locale} $X$, omitting the explicit mention of the base $B$.
A \emph{continuous map} of locales $f:X\to Y$ is a class function $f^-:B_Y\to X$ (note the reverse direction) satisfying:

\begin{itemize}
\item [1.] $\bigvee_{a\in B_Y}f^-(a)= \top$,
\item [2.] $f^-(a) \wedge f^-(b)= \bigvee \{f^-(c): c\in B_Y, c\leq  a, c\leq b\}$, for all $a,b\in B_Y$,
\item [3.] $f^-(a)\leq \bigvee_{b\in U} f^-(b)$, for all $a\in B_Y, U\in \Pow(B_Y)$ with $a\leq \bigvee U$.
\end{itemize}
\noindent By Replacement, one may assume that any such function is a set.
The (in general proper) class of these maps is denoted by $Hom(X,Y)$. Observe  that $Hom(X,Y)$ is in a one-to-one correspondence with the collection of frame homomorphisms from $Y$ to $X$ (i.e., class functions from $Y$ to $X$ preserving the frame structure).

A locale $X$ is \emph{compact} iff every covering of $\top$ by basic elements (i.e., every $U\in \Pow(B)$ such that $\top= \bigvee U$)  has a  finite  subcover. $X$ is \emph{regular} if, for all $a\in B$, $a= \bigvee \{b\in B: b\prec a\}$, where, for $x,y\in X$, $y\prec x\iff \top = x\vee y^*$, with $y^*= \bigvee \{c\in B: c\wedge y=\bot\}$  the \emph{pseudocomplement} of $y$.

Let $rc:B\to \Pow(B)$ be a function with the property that for all $b\in rc(a)$ a scale exists from $b$ to $a$, i.e., a family $\{c_p\}_{p\in \mathbb{I}}$ of elements of $X$, indexed on the rational unit interval $\mathbb{I}$,  satisfying: $c_0=b, c_1=a$ and, for $p< q$, $c_p\prec c_q$. $X$ is \emph{completely regular } if a function  $rc:B\to \Pow(B)$ of this kind is given with the property that for all $a\in B$, $a= \bigvee rc(a)$.  Over CZF (or HHA) plus the principle of dependent choice, a compact regular locale is completely regular.
In ZFC, compact (completely) regular locales and compact Hausdorff spaces define equivalent categories  \cite{J82}.

A locale $X$ is \emph{set-presented} \cite{A,CSSV} iff there is a function $C:B\to \Pow(\Pow(B))$, with $\Pow(\Pow(B))$
the class of subsets of the class $\Pow(B)$,
such that $$a\leq \bigvee U \iff \exists V\in C(a)\; V\subseteq U.$$
Observe that
this implies $a\leq \bigvee V$, for all $a\in B, V\in C(a)$.


For $x,y\in X$, one says that $y$ is \emph{way-below} $x$  iff for all $U\in \Pow(B)$, $x\leq \bigvee U$ implies
$y\leq \vee u$ for $u$ a finite subset of $U$. Due to the unbounded quantification over  $\Pow(B)$, the way-below relation is a class in CZF even when restricted to basic elements. However, if  $X$ is set-presented, for $a,b\in B$, one has $b$ \emph{way-below} $a$  iff $b<\!\!<a \equiv (\forall V\in C(a))(\exists v\in \Pow_{fin}(V)) b\leq \vee  v$. Since in CZF (or CZF$_{exp}$), the class $\Pow_{fin}(A)$ of  finite subsets of a given set $A$ is a set, the last formula is easily seen to be equivalent to a bounded formula, and thus defines a subset of $B\times B$.
A \emph{locally compact} locale  is a set-presented locale $X$ satisfying  $a= \bigvee \{b\in B: b<\!\!<a\}$, for all $a\in B$.
A key fact for what follows is that, in a locally compact locale, for $x,y\in X$, $y$ way-below $x$ implies $y$ way-below $\vee u$ and  $\vee u$ way-below  $x$, for $u$ a finite subset of $B$ \cite{CuSC}.

Finally, every compact regular locale is locally compact, with $<\!\!<=\prec$, and  is therefore  set-presented in CZF \cite{J82,A}.

The reader may consult \cite{J82} for the definition of the locale of the real numbers  $R$ and of its sublocale  $[0,1]$, the localic real unit interval.
Both locales are completely regular, $R$ is locally compact, $[0,1]$ is compact (hence both are set-presented in CZF).

The following result has an independent interest, in particular in connection with the theory of (rings of) continuous real-valued functions.

\begin{thm}[CZF]\label{smallHom}
If $X$ is a compact regular locale and $Y$ is a regular and set-presented locale, the class $Hom(X,Y)$ is a set.
\end{thm}
\noindent \textbf{Proof.} Note first that, for $X$ any locale and $Y$ a set-presented locale,
a mapping  $f^-:B_Y\to X$ satisfies conditions from 1 to 3 on continuous mappings iff it satisfies 1, 2 and

\begin{itemize}
\item [$3'$.] $f^-(a)\leq \bigvee_{b\in V} f^-(b)$, for all $a\in B_Y, V\in C(a)$.
\end{itemize}

\noindent
Given any continuous map $f\in Hom(X,Y)$, the associated morphism  $f^-:B_Y\to X$ is easily seen to preserve the well-inside relation $\prec$, in the sense that  $f^-(b) \prec  f^-(a)$, for all $a,b\in B_Y$ with $b\prec a$. As we assumed here  $Y$ to be regular, we also have $f^-(a)= \bigvee_{b\prec a} f^-(b)$.
Since in a compact regular locale the way-below relation coincides with the well-inside relation, for all $a,b\in B_Y$ with $b\prec a$,  there is  a finite subset $u$ of $B_X$ such that $f^-(b) <\!\!< \vee u <\!\!<  f^-(a)$. In particular, $f^-(b) \leq \vee u \leq  f^-(a)$.

Now let $W$ be the set $\{(a,b)\in B_Y \times B_Y: b\prec a\}$, and  consider the class $mv(\Pow_{fin}(B_{X})^W)$ of multivalued functions from $W$ to the set of finite subsets of $B_{X}$. By Fullness there is a set  $$K\subseteq mv(\Pow_{fin}(B_{X})^W)$$  such that for all  $R\in mv(\Pow_{fin}(B_{X})^W)$ there is $R^0\in K$ with $R^0\subseteq R$.
There is then a class function from $K$ to the class $X^{B_Y}$ of mappings from $B_Y$ to $X$,  defined by letting, for $F\in K$ and $a\in B_Y$, $f_F^-(a)= \bigvee(\bigcup_{b\prec a}\{\vee u:F((a,b),u)\})$. By Replacement, the range  $K'\subseteq X^{B_Y}$ of this function is a set.
We show that $Hom(X,Y)$ is a subset of this set. Indeed, let $D$ be the subclass of the set $K'$ defined by $f^-\in D$ iff $f^-$ satisfies conditions 1, 2 and $3'$ on continuous mappings. As $f^-$ is restricted to belong to the set $K'$, using Replacement these conditions can equivalently be expressed by a bounded formula, so that by Restricted Separation $D$ is a set.

Clearly, $D\subseteq Hom(X,Y)$.
Conversely, given $f\in Hom(X,Y)$, let the set $R_f\subseteq W\times \Pow_{fin}(B_X)$ be defined by:
$$R_f((a,b),u)\iff f^-(b) \leq \vee u \leq  f^-(a).$$
We have shown that for all $(a,b) \in W$ there is $u\in \Pow_{fin}(B_X)$ such that $R_f((a,b),u)$, so that $R_f\in mv(\Pow_{fin}(B_{X})^W)$. Then there is $R^0_f\in K\subseteq mv(\Pow_{fin}(B_{X})^W)$ such that $R^0_f\subseteq R_f$. Since
it also holds $(\forall (a,b) \in W)(\exists u \in \Pow_{fin}(B_{X})) R^0_f((a,b),u)$, by $f^-(a)= \bigvee_{b\prec a} f^-(b)$, one obtains  $$f^-(a)= \bigvee(\bigcup_{b\prec a}\{\vee u:R^0_f((a,b),u)\}),$$ for every $a\in B_Y$. Thus, $f^-= f^-_{R^0_f}\in K'$, and since $f^-$ satisfies conditions 1, 2 and $3'$, $f^-\in D$.

\bigskip

\noindent
Recall that
compact (completely) regular locales have the role in intuitionistic settings that compact Hausdorff spaces play in ordinary classical topology \cite{J82}. Let \textbf{KRLoc} (resp. \textbf{KCRLoc}) denote the full subcategory of the category \textbf{Loc} of locales whose objects are the compact regular (resp. compact completely regular) locales. By the above theorem we directly have:

\begin{cor}[CZF]
\textbf{KRLoc} and (hence) \textbf{KCRLoc} are locally small categories.
\end{cor}

\begin{cor}[CZF]\label{HOMXR}
If $X$ is a compact regular locale, then the classes $Hom(X,$ $ R)$ and $Hom(X,[0,1])$ are sets.
\end{cor}

\noindent \textbf{Remarks.}  $i.$ Note that the last corollary does not hold in CZF$_{exp}$: for $X=\Pow(\{1\})$,  $Hom(X,R)$ is isomorphic with the class of Dedekind reals which form a proper class in CZF$_{exp}$ \cite{LubarskyRathjen}.
However, Corollary \ref{HOMXR} can be proved in CZF$_{exp}$ plus the principle of countable  choice AC$_\omega$.

$ii.$ In \cite{CuSC} I proved that in stronger systems, as CZF+uREA+DC, or CTT, the class $Hom(X,Y)$ is a set whenever $X$ is locally compact and $Y$ is set-presented and regular. Using the  type-theoretic axiom of choice and regular universes, Palmgren \cite{Pa} generalized this result in CTT by weakening $X$ locally compact to $X$ set-presented.

$iii.$ Over CZF+REA, the concept of set-presented locale and of \emph{inductively generated} formal space \cite{CSSV} are equivalent, cf.  \cite{A}. Although this is no more the case over CZF,  one may prove with essentially the same argument that Theorem \ref{smallHom} also holds  if one replaces (locales with formal spaces and) `set-presented' with `inductively generated'.

$iv.$ Finally, we note that the proof of Theorem \ref{smallHom} can easily be seen to apply in the more general case where  $X$ is any locally compact locale and $Hom(X,Y)$ is replaced by the class of mapping $f$ in $Hom(X,Y)$ such that $f^-$ sends the well-inside relation on a base of $Y$ in the way-below relation (i.e., $f^-(b)<\!\!<  f^-(a)$ whenever $b\prec a$, $a,b\in B_Y$). A further generalization of a more formalistic nature is implicit in the proof of Theorem \ref{smallHom}, but  will not be spelt out here.

\paragraph{Existence of Stone-\v{C}ech compactification.}\setcounter{subparagraph}{0}
The (generalized) Sto\-ne-\v{C}ech compactification  of a space or locale $X$ is  its compact completely regular reflection, i.e., it is a continuous map  $$\eta:X\to \beta X,$$ with $\beta X$ compact and completely regular, which  satisfies  the following universal property: for every continuous map $f:X\to Y$ to a compact completely regular codomain $Y$, a unique map $f^\beta:\beta X \to Y$ exists such that $f^\beta \circ \eta= f$. Moreover, $\eta$ is a (dense) embedding precisely when $X$ is completely regular.

In \cite{CuSC} one finds the proof in CTT and CZF + REA that, on the assumption that for every compact completely regular locale $Y$ the class $Hom(Y,[0,1])$ is a set,  if $X$ is such that $Hom(X,[0,1])$ is a set, the Stone-\v{C}ech compactification of $X$ exists \cite[Corollary 6.2]{CuSC} (see also \cite{CuSCAlexandroff}).

By Theorem \ref{smallHom} we may then conclude that if $Hom(X,[0,1])$ is a set then $\beta X$ exists in CZF + REA, in particular without any intervention of a choice principle.
As, moreover, the universal property of $\beta$ directly yields a bijection $Hom(X,[0,1])\cong Hom(\beta X,[0,1])$, and since the latter is a set by Theorem \ref{smallHom} when $\beta X$ exists, in summary we have:

\begin{thm}[CZF+REA]
The Stone-\v{C}ech compactification $\beta X$ of a locale $X$ exists  if, and only if, $Hom(X,[0,1])$ is a set.
\end{thm}

The universal property of $\beta$ yields more generally a bijection $Hom(X,Y)\cong Hom(\beta X,Y)$, for every compact completely regular $Y$, and, again by Theorem  \ref{smallHom}, one has that $Hom(\beta X,Y)$ is a set. Thus:

\begin{cor}[CZF+REA]
Given any locale $X$, if $\beta X$ exists, $Hom(X,Y)$ is a set for every compact completely regular locale $Y$. Therefore, for every locale $X$,  $Hom(X,Y)$ is a set for every compact completely regular locale $Y$ if and only if $Hom(X,[0,1])$ is a set.
\end{cor}

\noindent In \cite{CuSC}, in the context of CZF+uREA+DC,  these results similarly followed  by the results recalled in  Remark $ii.$ of the previous section.

\medskip

We now turn to the proof that there are locales $X$ of which $\beta X$ does not exist constructively. A well-known classical identification is that of the subsets of a set $S$ with the mappings of $S$ in the two-element set $\{0,1\}$. We may regard this identification, that fails intuitionistically, as a special case of a bijection that exists between the frame of opens $\Omega(X)$ of a Boolean space $X$ (i.e., a space $X$ such that $\Omega(X)$ is Boolean), and the continuous functions from $X$ to $\{0,1\}$. The next lemma shows that  this classically more general fact does constructively carry over to the localic setting.
 Recall that a locale is \emph{Boolean} if it is a (complete) Boolean algebra.

\begin{lem}[CZF$_{exp}$]\label{01classifies}
The compact completely regular locale of opens of the discrete space $\{0,1\}$ classifies the opens of any Boolean locale $X$. I.e., a bijection exists between the class of elements of $X$ and $Hom(X, \Pow(\{0,1\}))$.
\end{lem}

\noindent \textbf{Proof.} To $a\in X$ one associates the map $f^-_a:\{\{0\},\{1\}\}\to X$, defined by $f^-_a(\{0\})=a$, $f^-_a(\{1\})=a^*$. We leave to the reader to check that this defines a continuous map.
Conversely,  $f:X\to \Pow(\{0,1\})$ defines the open $a_f= f^-(\{0\})$.

\medskip

\noindent Note that
$\Pow(\{0,1\})$ is not Boolean intuitionistically. To prove that there are locales of which Stone-\v{C}ech compactification does not exist, we shall need the result that the class  of continuous mappings from a compact completely regular locale to $\Pow(\{0,1\})$ is a set. In CZF, this follows by Theorem \ref{smallHom}, as $\Pow(\{0,1\})$ is set-presented and regular. However,  this also holds in CZF$_{exp}$.

\begin{lem}[CZF$_{exp}$]\label{homcsp01}
Let $X$ be any compact locale. Then the class $Hom(X,$ $ \Pow(\{0,1\}))$ is a set.
\end{lem}

\noindent \textbf{Proof.} Note first that, since $\Pow(\{0,1\})$ is set-presented, $f\in Hom(X,$ $ \Pow(\{0,1\}))$ iff $f^-:\{\{0\},\{1\}\}\to X$ satisfies conditions  1, 2 and $3'$ on continuous mappings (cf. the proof of Theorem \ref{smallHom}). By Replacement, the class $\bar B_X= range(\vee)$, $\vee:\Pow_{fin}(B_{X})\to X$, is a set, so that by Exponentiation also the class $\bar B_X ^{\{\{0\},\{1\}\}}$ is  a set. We show that $Hom(X, \Pow(\{0,1\}))$ coincides with the subclass $D$ of $\bar B_X ^{\{\{0\},\{1\}\}}$ given by the maps in $\bar B_X ^{\{\{0\},\{1\}\}}$ satisfying conditions 1, 2 and $3'$. As $\bar B_X ^{\{\{0\},\{1\}\}}$ is a set, exploiting the Replacement Scheme these conditions can equivalently be  expressed by a bounded formula, so that $D=Hom(X, \Pow(\{0,1\}))$ is a set by Restricted Separation.

Let then $f\in Hom(X, \Pow(\{0,1\}))$ be a continuous map. One must have $\top= f^-(\{0\})\vee f^-(\{1\})$, and $\bot= f^-(\{0\})\wedge f^-(\{1\})$. As $X$ is compact, there is a finite subset $v$ of $\{b\in B_X:b\leq f^-(\{0\})\}\cup  \{b\in B_X:b\leq f^-(\{1\})\}$ such that $\top=\vee v$. There are thus finite subsets $v_0,v_1$ such that $v_0\cup v_1 = v$, $v_0\subseteq \{b\in B_X:b\leq f^-(\{0\})\}$, and $v_1\subseteq \{b\in B_X:b\leq f^-(\{1\})\}$. Hence, $\top=\vee (v_0\cup v_1)$, and $\bot = (\vee v_0) \wedge (\vee v_1)$. This gives $f^-(\{0\})=\vee v_0$ and $f^-(\{1\})=\vee v_1$. Thus,  $f^-$ is in fact a map from $\{\{0\},\{1\}\}$ to $\bar B_X$ satisfying conditions  1, 2 and $3'$, i.e., $f^- \in D$. The converse inclusion is trivial.

\medskip
\noindent These two lemmas, together with Lemma \ref{Noframeset}, give us:

\begin{thm}\label{noSCBoole}
The Stone-\v{C}ech compactification  of a non-degenerate Boolean locale $X$ cannot be defined in CZF (+REA+PA+Sep). Moreover, $Hom(X, $ $\Pow(\{0, 1\}))$, and $Hom(X, [0,1])$, $Hom(X, R)$, are proper classes in this setting.
\end{thm}

\noindent \textbf{Proof.} By Lemma \ref{homcsp01} (or Theorem \ref{smallHom}), if $\beta X$ exists, $Hom(\beta X,\Pow(\{0,1\}))$ is a set in CZF (and a fortiori in every stronger system).  Moreover, by the universal property of $\beta$, $Hom(X,\Pow(\{0,1\})\cong Hom(\beta X,\Pow(\{0,1\})$. Thus  $Hom(X,\Pow(\{0,1\})$ is a set, so that, by Lemma \ref{01classifies}, $X$ is a set in CZF$_{exp}$.
However, by Lemma  \ref{Noframeset}, no non-degenerate locale can be proved to have a set of elements in CZF (+REA+PA+Sep).

\medskip

\noindent  Boolean locales abound in nature also constructively: given any (non-degene\-rate) locale $X$, the set-generated class-frame of its regular elements \cite{FS79} is a (non-degenerate) Boolean locale. Note also that a Boolean locale is (completely) regular.

\medskip

\noindent \textbf{Remarks.}
$i.$ As, e.g. in IZF, which is obtained from CZF by replacing Restricted Separation with Separation and adding Powerset, Stone-\v{C}ech compactification of every locale can be defined, Theorem \ref{noSCBoole} is an independence result.

$ii.$ In the statements of the Main Lemma and of Theorem \ref{noSCBoole} one may of course replace CZF+REA+PA+Sep with any extension of CZF (or even of CZF$_{exp}$) that is compatible with GUP.
Similarly, the above results can be shown to hold (mutatis mutandis) also with respect to every extension of  CTT compatible with (the  type-theoretic version of) GUP.

$iii.$ The given proof of Theorem \ref{noSCBoole} is entirely self-contained (it does not depend on results in \cite{CuSC});  a  corresponding result is analogously seen to hold also for the  compact zero-dimensional reflection  of a Boolean locale, as a zero-dimensional locale is regular, and as $\Pow(\{0,1\})$ is (compact and) zero-dimensional.

$iv.$ By contrast with Theorem \ref{noSCBoole}, the `approximation' to  Stone-\v{C}ech compactification introduced in \cite{CuSC} exists for every locale $L$ (and every given set-indexed family of continuous maps of the appropriate type, see \cite{CuSC}).

\paragraph{Conclusion.} If one agrees in considering a necessary condition in order for an argument to be defined constructive that it may be formulated within an extension of  CZF or CTT compatible with the  form of the uniformity principle we are considering,
Theorem \ref{noSCBoole} is read as saying that the Stone-\v{C}ech compactification  of a non-degenerate Boolean locale $X$ does not exist constructively. This goes very much against what
holds in the topos-theoretic context: in any topos, the Stone-\v{C}ech compactification of a Boolean locale $X$ is simply given by the lattice of ideals on $X$.

Note that the given necessary condition for constructivity is by no means sufficient: the theory CZF+Sep, where Sep is impredicative unbounded separation, has the same proof-theoretic strength of second-order Heyting arithmetic \cite{Lubarsky2006}, and is however consistent with the generalized uniformity principle. We find it remarkable, and somewhat surprising, that, due to this fact, Lemma \ref{Noframeset} and Theorem \ref{noSCBoole} also hold with respect to this theory.

{\small

\bibliographystyle{plain}

}
\end{document}